%% file: 1995-066.tex
\input amstex
\input epsf
\documentstyle{amsppt}
	\magnification=1200
	\rightheadtext{Random walks of circle packings}

\input randomw_defs.tex

\def\today{\number\day\space\ifcase\month\or January \or
February \or March \or April \or May \or June \or July \or
August \or September \or October \or November \or
December\fi\space\number\year}


\topmatter
\title 
Random walks of circle packings
\endtitle
\author
Tomasz Dubejko
\endauthor
\address
Mathematical Sciences Research Institute, Berkeley, CA 94720\endaddress
\email
tdubejko$\@$msri.org\endemail
\thanks
Research at MSRI is supported in part by grant no.DMS-9022140.\endthanks
\abstract
A notion of random walks for circle packings is introduced.
The geometry behind this notion is discussed, together with some applications.
In particular, we obtain a short proof of a result regarding the type problem for circle packings, which shows that the type of a circle packing is closely related to the type of its tangency graph. 
\endabstract
\endtopmatter

\document

\heading 1.~Introduction \endheading

Circle packings have been introduced by Thurston to build finite Riemann mappings and to conjecture the Finite Riemann Mapping Theorem ([T1], [T2]).
They have been used to construct discrete analytic function theory ([BSt1], [D1], [D2], [HR], [RS], [St1]).
Interesting connections between \cp s and graph-theoretical aspects of their tangency patterns have been investigated in [St1], and recently in [Mc], [BeSc], [D3], and [HSc2].

Here we are interested in developing geometric notion of random walks and harmonic functions for \cp s.
This notion will be first described for \cp s with mutually disjoint interiors, i.e. univalent packings, and then will be extended to branched and overlapping-angles \cp s.
Our construction will yield explicit harmonic functions for random walks of packings, which are also Dirichlet-finite when associated packings have compact carriers.
It will also follow that the real and imaginary parts of finite Riemann mappings of [T1] are harmonic for underlying random walks.
Moreover, we obtain short (and elementary, as they can be) proofs of the following.

\proclaim{Theorem~1.1}
Let \mcompk\ be a triangulation of an open disc in the complex plain $\bold C$.
If the simple random walk on the 1-skeleton \medgek\ of \mcompk\ is recurrent then there exists a univalent \cp\ for \mcompk\ whose carrier is $\bold C$.
Conversely, if \mcompk\ is of bounded degree and there exists a univalent \cp\ for \mcompk\ whose carrier is the unit disc $\bold D$ then the simple random walk on \medgek\ is transient.
\endproclaim

\proclaim{Theorem~1.2}
Let \mcompk\ be a triangulation of an open disc.
If \medgek\ is of bounded degree and transient then there exist non-constant, Dirichlet-finite, harmonic functions for the simple random walk on \medgek.
\endproclaim

Several remarks regarding the above theorems should be made.
It was shown ([BSt1], [HSc1]) that given \mcompk\ as in Theorem~1.1 there exists a univalent \cp\ for \mcompk\ whose carrier is $\bold C$ or $\bold D$, but not both.
The problem of deciding whether the carrier is $\bold C$ or $\bold D$ has become to be known as the type problem for \cp s (see [BSt1]), by analogy to the classical type problem for graphs.
The first answer regarding this problem was given in [Mc] for \mcompk\ of bounded degree.
Theorem~1.1 subsume the result of [Mc], and has been very recently established in [HSc2] using interesting but involved combinatorial extremal length techniques.
We give here simple and independent proof of this result as an example of an application of our new notions.

The advantage of working with random walks of packings, rather then with simple random walks of their tangency graphs, is that there exist non-trivial harmonic functions associated with them which have a geometric interpretation.
Using geometric properties of these functions and the fact [Me] that the existence of Dirichlet-finite harmonic functions is preserved by rough isometries of networks, Theorem~1.2 easily follows from Theorem~1.1.
We note that Theorem~1.2 was originally proved in [De] using cohomology theory and also follows from a new result of [BeSc].
We also observe that Theorem~1.2 and [So2]/[Me] imply that there exist non-constant, Dirichlet-finite, harmonic functions for graphs of bounded degree which are roughly isometric to the 1-skeleton of a triangulation of Theorem~1.2;
this includes, for example, planar trees of bounded degree and planar cell decompositions of bounded degree with a uniform bound on the number of edges in each cell.

\remark{Acknowledgments}
The author would like to thank K. Stephenson for interesting discussions and P. M. Soardi for the reference [Me].
Figure~2 was created with the help of Stephenson's software {\tt CirclePack}.
\endremark

\heading 2.~Preliminaries \endheading

Since there are number of publications on the subject of \cp s, we skip formal definitions.
Our terminology will be consistent with that of [D1] and [DSt] (see also [BSt1]) for \cp s without overlapping-angles.
The reader is refered to [BoSt] and [T2] for information about \cp s with angles of overlaps.

We will now set up our notation.
If \mcompk\ is a triangulation, finite or infinite, with or without boundary, of a disc in $\bold C$ then we will write \mvertk\ and \medgek\ for sets of vertices and edges of \mcompk, respectively.
If \mmp P is a \cp\ for \mcompk\ then \mcircle Pv will denote the circle in \mmp P corresponding to a vertex $v\in \vertk$, \mr Pv will be the radius of \mcircle Pv, and $\smap P: \compk \to \carr P$ will be the simplicial isomorphism of \mcompk\ with the carrier \mcarr P of \mmp P determined by: $\smap P(v):=\text{the center of \mcircle Pv}$, $v\in \vertk$.

As far as the subject of random walks is concerned, the survey article [W] is more than a sufficient source of information for our purposes (see also [So1]).

If \mcompk\ is as above then $(\compk, p)$ will denote a random walk on the 1-skeleton of \mcompk\ (i.e., on the graph \medgek), where $p: \vertk \times \vertk \to [0,1]$ is the transition probability function.
To simplify the notation, if $u,v\in \vertk$ are adjacent then we will write $u\sim v$, and we will denote their common edge by $uv$.
In the sequence we will be working with random walks given by conductances.
If $c:\vertk \times \vertk \to (0,\infty)$ is such that
\roster
\item
$c(u,v)>0$ iff $u\sim v$,
\item
$c(u,v)=c(v,u)$,
\item
$c(v,v)=0$ for all $v\in \vertk$,
\endroster
then the random walk $(\compk, p)$ given by $p(u,v):=c(u,v)/(\sum_{w\sim u}c(u,w))$ will be called {\it revesible} with the {\it edge conductance} between $u$ and $v$, $u\sim v$, equal to $c(u,v)$.
The {\it simple} random walk on \medgek\ is a random walk with all edge conductances equal to 1.
A real-valued function $f:\vertk \to \bold R$ is said to be {\it harmonic} for $(\compk, p)$ if 
$$
f(v)=\sum_{u\sim v}p(v,u) f(u) \quad \text{for every $v\in \intk$},
$$
where \mintk\ is the set of interior vertices of \mcompk.
A function $f:\vertk \to \bold R$ is said to be {\it Dirichlet-finite} for a reversible random walk $(\compk, p)$ with edge conductances $c(u,v)$ if
$$
\sum_{uv\in \edgek} \bigl( f(u) - f(v) \bigr)^2 c(u,v) <\infty .
$$
A random walk $(\compk, p)$ is called {\it recurrent} if probability of visiting any vertex infinitely many times is 1; otherwise is called {\it transient}.
Finally, a graph \medgek\ is called {\it recurrent} ({\it transient}) if the simple random walk on \medgek\ is recurrent (transient).

Two important results linking the existence of harmonic functions with the type of a random walk are summarized in the following.

\proclaim{Theorem~2.1}{\rm [W, Thm.~4.1, Cor.~4.15]}
Let \mcompk\ be a triangulation of an open disc in $\bold C$.
\roster
\item
If $(\compk, p)$ is recurrent then there are no non-constant, non-negative harmonic functions for $(\compk, p)$.
\item
If $(\compk, p_1)$ and  $(\compk, p_2)$ are two reversible random walks with edge conductences $c_1$ and $c_2$, respectively, such that $c_1(u,v)\ge c_2(u,v)$ for all $uv\in \edgek$ then the recurrence of $(\compk, p_1)$ implies the recurrence of $(\compk, p_2)$.
\endroster
\endproclaim

We also need to recall the notion of rough isometry.
Two metric spaces $(X, d_X)$ and $(Y, d_Y)$ are said to be {\it roughly isometric} if there exist a map $\psi: X\to Y$ and constants $A>0$ and $B\ge 0$ such that
$$
A^{-1}d_X(x,x') - B\le d_Y(\psi(x), \psi(x'))\le Ad_X(x,x') + B \quad \text{for all $x,x'\in X$},
$$
and
$$
d_Y(y, \psi(X))\le B \quad \text{for every $y\in Y$}.
$$

For a graph $\Bbb X$ and a reversible random walk $(\Bbb X,p)$ on $\Bbb X$ with edge conductances $c(\cdot ,\cdot)$, we define the associated metric space $(\Bbb X, d_p)$ as follows.
First we define the length of an edge path $\gamma$ in $\Bbb X$ by
$$
\ell(\gamma):= \sum_{xx'\in \gamma} \frac{1}{c(x,x')}.
$$
Then the metric $d_p$ on $\Bbb X$ is given by
$$
d_p(x,x')=\inf_{\gamma} \{ \ell(\gamma):\ \text{$\gamma$ is a path in $\Bbb X$ with end points $x$ and $x'$} \}.
$$

We will say that two reversible random walks $(\Bbb X,p)$ and $(\Bbb Y,q)$ are roughly isometric if $(\Bbb X,d_p)$ and $(\Bbb Y,d_q)$ are roughly isometric.
We finish preliminaries with the following result that has been shown in [Me] (see also [So2]).

\proclaim{Theorem~2.2}{\rm [Me, Thm.~3.6]}
Let $\Bbb X$ and $\Bbb Y$ be two graphs of bounded degree.
Suppose $(\Bbb X,p)$ and $(\Bbb Y,q)$ are roughly isometric reversible random walks such that $c_p(x,x')<\kappa$ and $c_q(y,y')<\kappa$ for some constant $\kappa$ and all edges $xx'$ and $yy'$, where $c_p$ and $c_q$ are edge conductances of $p$ and $q$, respectively.
Then there exists a non-constant, Dirichlet-finite, hermonic function for $(\Bbb X,p)$ iff there exists one for $(\Bbb Y,q)$.
\endproclaim

\heading 3.~Random walks of univalent packings \endheading

In this section we will be working only with univalent \cp s, i.e. \cp s whose circles have mutually disjoint interiors.
Suppose \mmp P is a univalent packing for \mcompk, and $v\in \intk$.
If $u$ is a neighbor of $v$ then there are exactly two vertices $w'$ and $w''$ that are adjacent to both $u$ and $v$.
Let $r_{\p P}'(uv)$ and $r_{\p P}''(uv)$ be radii of two orthogonal circles to triples $\circle Pu, \circle Pv, C_{\p P}(w')$ and $\circle Pu, \circle Pv, C_{\p P}(w'')$, respectively (Fig.~1).

\midinsert
\epsfysize=8.5truecm
\centerline{\epsffile[219 310 390 480]{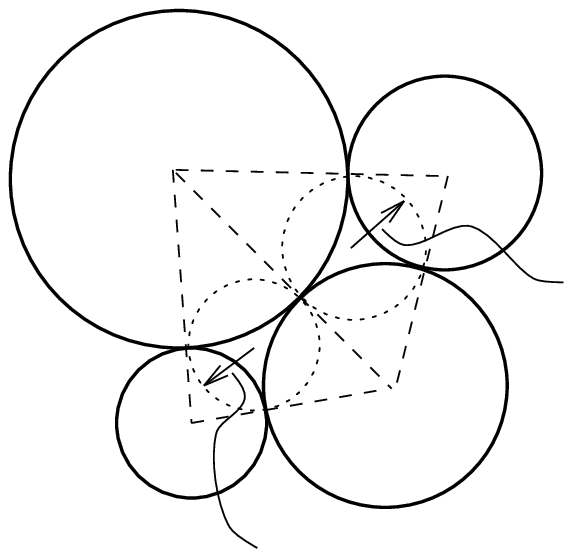}}
\vskip-8.0truecm\hskip5.0truecm $\circle Pu$ 
\vskip0.0truecm\hskip10.5truecm $C_{\p P}(w'')$
\vskip2.3truecm\hskip11.3truecm $r_{\p P}''(uv)$
\vskip2.3truecm\hskip3.5truecm $C_{\p P}(w')$ \hskip3.0truecm $\circle Pv$
\vskip0.5truecm\hskip6.9truecm $r_{\p P}''(uv)$
\captionwidth{25pc}
\botcaption{Figure~1} Orthogonal circles.
\endcaption
\endinsert

We define the edge conductance between $u$ and $v$ induced by \mmp P by
$$
c_{\p P}(u,v):=\frac{r_{\p P}'(uv)+r_{\p P}''(uv)}{\r Pu + \r Pv}. 
\tag"(*)"
$$
We can now introduce

\definition{Definition~3.1}
If \mmp P is a univalent \cp\ for \mcompk\ then the random walk $(\compk, p_{\p P})$ induced by \mmp P is defined as follows:
$$
p_{\p P}(u,v)=\cases \frac{c_{\p P}(u,v)}{\sum_{w\sim u}c_{\p P}(u,w)}, &\text{if $v\in \intk$, $u\sim v$, and $u\ne v$}\\
	0, &\text{if $v\in \intk$ and $u=v$, or $v\notin \intk$ and $u\ne v$}\\
	1, &\text{if $v\notin \intk$ and $u=v$}.\endcases
$$
\enddefinition

\remark{Remark~3.2}
\roster
\item"1)"
If \mcompk\ has no boundary then $(\compk, p_{\p P})$ is reversible with edge conductances $c_{\p P}(u,v)$.
If \mcompk\ has boundary then the boundary is an absorbing barrier for $(\compk, p_{\p P})$.
\item"2)"
From Herron's formula it follows that $c_{\p P}(u,v)$ can be written explicitly in terms of radii of circles of \mmp P,
$$
c_{\p P}(u,v)=\tfrac{\sqrt{\r Pu \r Pv}}{\r Pu + \r Pv}\left( \sqrt{ \tfrac{r_{\p P}(w')}{\r Pu + \r Pv + r_{\p P}(w')} } + 
\sqrt{ \tfrac{r_{\p P}(w'')}{\r Pu + \r Pv + r_{\p P}(w'')} }\right),
$$
where $w'$ and $w''$ are vertices adjacent to both $u$ and $v$.
\item"3)"
If \mcompk\ has no boundary and is of bounded degree, and \mmp P is a univalent packing for \mcompk\ then the Ring Lemma [RS] implies that there is a constant $\kappa$, depending only on the degree of \mcompk, such that $c_{\p P}(u,v)\in (\kappa^{-1}, \kappa)$ for all edges $uv\in \edgek$.
In particular, by Theorem~2.1.(2), $(\compk, p_{\p P})$ is recurrent iff \medgek\ is recurrent.
Moreover, $(\compk, p_{\p P})$ and the simple random walk on \medgek\ are roughly isometric.
\endroster
\endremark

We will now show that with every random walk induced be a \cp\ there are associated two non-trivial harmonic functions.

\proclaim{Lemma~3.3}
Suppose \mmp P is a univalent \cp\ for \mcompk\ and $(\compk, p_{\p P})$ is the induced random walk on \medgek.
Let $\Re_{\p P}:\vertk \to \bold R$ and $\Im_{\p P}:\vertk \to \bold R$ be defined by $\Re_{\p P}(v):=\text{{\rm Re}}\bigl(\smap P(v)\bigr)$ and $\Im_{\p P}(v):=\text{{\rm Im}}\bigl(\smap P(v)\bigr)$, respectively, where $\smap P:\compk \to \carr P\subset \bold C$ is the associated simplicial isomorphism, and {\rm Re} and {\rm Im} denote real and imaginary parts of a complex number, respectively.
Then functions $\Re_{\p P}$ and $\Im_{\p P}$ are harmonic for $(\compk, p_{\p P})$.
\endproclaim

Before we give a proof of the lemma we need the following proposition.

\proclaim{Proposition~3.4}
Let $z_0,\dots,z_m \in \bold C$ and let $\gamma = z_0z_1\cup z_1z_2\cup \dots \cup z_mz_0$ be a closed, piecewise-linear, oriented path in $\bold C$ built out of oriented segments $z_jz_{j+1}$ ($j$~\text{{\rm mod}}~$(m+1)$).
Denote by $\eta_j$ the unit normal vector to $z_jz_{j+1}$.
Then
$$
\sum_{j=0}^m |z_{j+1}-z_j|\eta_j =0,
$$
where $|z|$ is the modulus of $z$.
\endproclaim

\demo{Proof}
Write $z_{j+1}-z_j = |z_{j+1}-z_j|e^{i\varphi_j}$, where $\varphi_j\in [0,2\pi)$. 
Then $\eta_j=e^{i(\varphi_j+\pi/2)}$.
Thus
$$
0=e^{i\pi/2}\bigl( \sum_{j=0}^m (z_{j+1}-z_j) \bigr) = \sum_{j=0}^m e^{i\pi/2}|z_{j+1}-z_j|e^{i\varphi_j} = \sum_{j=0}^m |z_{j+1}-z_j|\eta_j.\qed
$$
\enddemo

\demo{Proof of Lemma~3.3}
It is sufficient to show that 
$$
\sum_{u\sim v} c_{\p P}(v,u)\bigl( \smap P(u) - \smap P(v) \bigr)=0 \quad \text{for every $v\in \intk$}.
$$
For $v\in \intk$ let $v^0,\dots,v^m$ be its heighboring vertices listed in the positive orientation.
Let $z_j$ be the center of the orthogonal circle to the triple of circles $C_{\p P}(v)$, $C_{\p P}(v^j)$, and $C_{\p P}(v^{j+1})$.
From (*) and Proposition~3.4 it follows
$$
\align
&\sum_{u\sim v} c_{\p P}(v,u)\bigl( \smap P(u) - \smap P(v) \bigr)= 
\sum_{u\sim v} c_{\p P}(v,v^j)\bigl( \smap P(v^j) - \smap P(v) \bigr)\\
&=\sum_{u\sim v} \frac{|z_j-z_{j-1}|}{r_{\p P}(v^j)+r_{\p P}(v)}\bigl( \smap P(v^j) - \smap P(v) \bigr)
=\sum_{u\sim v} |z_j-z_{j-1}|\frac{\smap P(v^j) - \smap P(v)}{r_{\p P}(v^j)+r_{\p P}(v)}= 0.
\endalign
$$
\qed
\enddemo

Suppose that \mmp P and \mmp Q are two univalent \cp s for \mcompk.
Then the \cp\ map $S_{\p P,\p Q}:=\smap Q\circ \smap P^{-1}:\carr P \to \carr Q$ maps centers of circles in \mmp P to centers of the corresponding circles in \mmp Q, and can be regarded as a discrete analog of an analytic function (see [BSt1], [D1], [RS], [T1]).
In particular, functions $\Re_{\p Q}$ and $\Im_{\p Q}$ can be viewed as a harmonic function and its conjugate for $(\compk, p_{\p Q})$, respectively.

In [RS] (see also [St1], [HR]) it was proved that finite Riemann mappings (which are \cp\ maps) of [T1] can be used to approximate the classical Riemann mapping.
If $\Omega$ is a Jordan domain in $\bold C$, $A\subset \Omega$ is compact, and $\epsilon>0$ then for every sufficiently small mesh regular hexagonal \cp\ \mmp P filling in $\Omega$ there exists combinatorially equivalent (i.e., with the same tangency graph) \cp\ \mmp Q filling in the unit disc $\bold D$ such that $|S_{\p P,\p Q}(z) - \tau(z)|<\epsilon$ for all $z\in A$, where $\tau$ is the Riemann mapping from $\Omega$ onto $\bold D$ (Fig.~2).
Using the fact that the ratio function $S_{\p P,\p Q}^{\#}:\carr P \to \bold R$, $S_{\p P,\p Q}^{\#}(v):=r_{\p Q}(v)/r_{\p P}(v)$, associated with $S_{\p P,\p Q}$ approximates $|\tau'|$ in $A$, it is not hard to show that the harmonic measure of $\partial A$ is approximated by escape (hitting) probabilities of $(\compk, p_{\p Q})|_{\compk_A}$, where \mcompk\ is the tangency graph of \mmp P (and \mmp Q), $\compk_A$ is the maximal simply connected 2-subcomplex of \mcompk\ such that $\{C_{\p P}(v)\}_{v\in \vertk_A}\subset A$, and $(\compk, p_{\p Q})|_{\compk_A}$ is the random walk on $\compk_A$ obtained by restricting $(\compk, p_{\p Q})$ to $\compk_A$ with $\partial \compk_A$ as absorbing barrier.
It is our conjecture that the harmonic measure of $\partial \Omega$ is approximated by escape probabilities of $(\compk, p_{\p Q})$.
However, this problem seems to be rather difficult to settle because up-to-date results show that $\tau$ and $|\tau'|$ can be approximated uniformly only on compacta of $\Omega$ by \cp\ maps.

\midinsert
\centerline{\epsfysize=6.5truecm\epsffile[90 155 500 610]{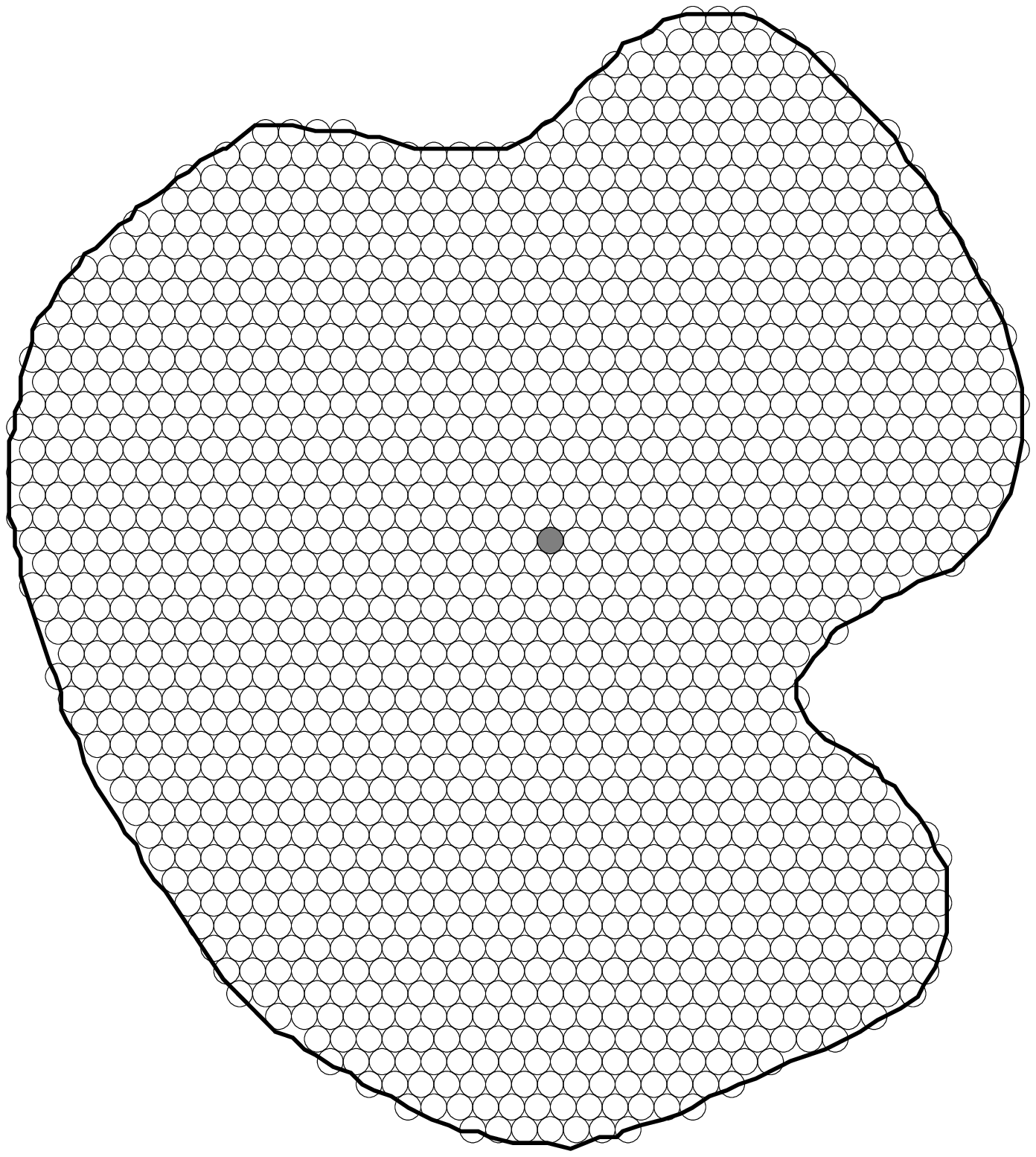} \hskip0.5truecm  \hskip2.5truecm \epsfysize=6.2truecm\epsffile[70 160 540 635]{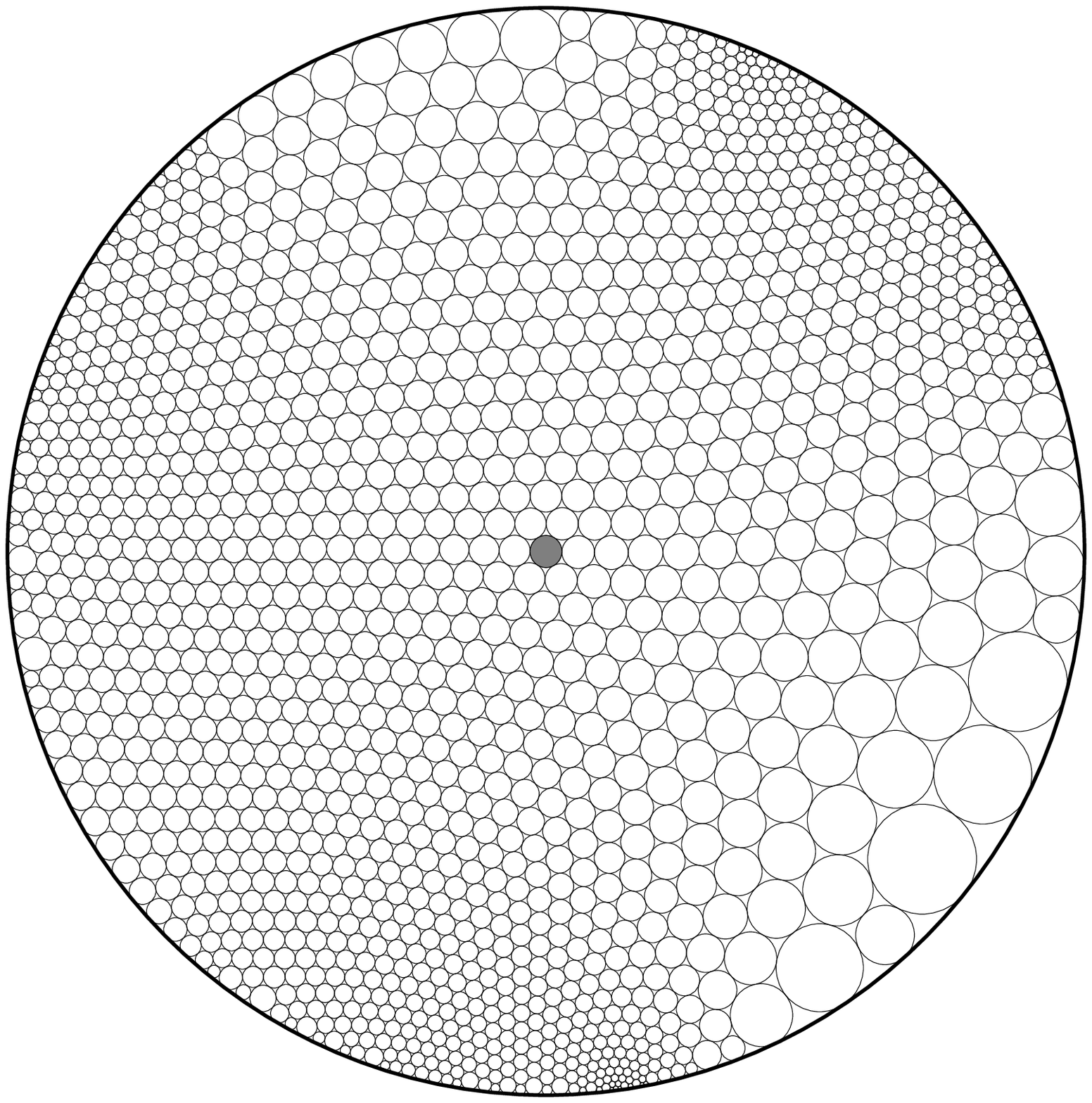}}
\vskip-3.5truecm\hskip6.0truecm \epsfxsize=1.5truecm\epsffile[255 385 355 405]{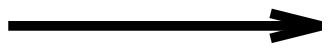}
\vskip2.5truecm\hskip4.2truecm $\Omega$ \hskip9.0truecm $\bold D$
\captionwidth{25pc}
\botcaption{Figure~2} Finite Riemann mapping.
\endcaption
\endinsert

We will now give proofs of Theorems~1.1 and 1.2.

\demo{Proof of Theorem~1.1}
Suppose \medgek\ is recurrent.
From [HSc1, Thm.~0.3] there exists either a univalent \cp\ for \mcompk\ with carrier equal to $\bold C$ or a univalent \cp\ with carrier equal to $\bold D$, but not both.
Denote this packing by \mmp P.
From Remark~3.2.2) it follows that
$$
c_{\p P}(u,v)\le 2\frac{\sqrt{r_{\p P}(u)r_{\p P}(v)}}{r_{\p P}(u)+r_{\p P}(v)} \le 1.
$$
Since edge conductances of the simple random walk on \medgek\ are all equal to 1, Theorem~2.1.(2) implies that $(\compk, p_{\p P})$ is recurrent.
Now, if \mcarr P were equal to $\bold D$ then, by Lemma~3.3, there would exist a bounded, non-constant, harmonic function for $(\compk, p_{\p P})$.
This, however, is impossible by Theorem~2.1.(1).
Hence $\carr P =\bold C$.

Suppose now that \mcompk\ is of bounded degree and there exists a univalent \cp\ \mmp P for \mcompk\ such that $\carr P =\bold D$.
To show that \medgek\ is transient it is sufficient, by Remark~3.2.3), to show that $(\compk, p_{\p P})$ is transient.
This follows immediately from Lemma~3.3 and Theorem~2.1.(1).
\qed
\enddemo

\remark{Remark~3.5}
The result in Theorem~1.1 is best possible.
It is easy to construct (see [BeSc]) a \cp\ with carrier equal to $\bold C$ whose tangency graph is transient (and, of course, of unbounded degree).
\endremark

\demo{Proof of Theorem~1.2}
From Theorem~1.1 and [HSc1, Thm.~0.3] there exists a univalent \cp\ \mmp P for \mcompk\ with $\carr P =\bold D$.
Let $\Re_{\p P}$ be as in Lemma~3.3.
Let $M$ be the degree of \medgek.
Then, by Remark~3.2.2),
$$
\align
&\sum_{uv\in \edgek}\bigl( \Re_{\p P}(u) - \Re_{\p P}(v) \bigr)^2 c_{\p P}(u,v)
\le \kappa\sum_{uv\in \edgek}\bigl( \Re_{\p P}(u) - \Re_{\p P}(v) \bigr)^2\\
&\le \kappa\sum_{uv\in \edgek}\bigl( r_{\p P}(u) + r_{\p P}(v) \bigr)^2
\le \kappa\sum_{uv\in \edgek} 2\bigl( r_{\p P}^2(u) + r_{\p P}^2(v) \bigr)
\le 2\kappa \sum_{v\in \vertk} M r_{\p P}^2(v)
\le 2\kappa M.
\endalign
$$
Hence, $\Re_{\p P}$ is a non-constant, Dirichlet-finite, harmonic function for $(\compk, p_{\p P})$.
By Theorem~2.2 and Remark~3.2.3), it follows that \mcompk\ has non-constant, Dirichlet-finite, harmonic functions.
\qed
\enddemo

From the last proof we obtain the following corollary.

\proclaim{Corollary~3.6}
Suppose \mcompk\ is a triangulation of an open disc in $\bold C$.
If \medgek\ is of bounded degree and transient, there exists a \cp\ \mmp P for \mcompk\ with \mcarr P compact.
Moreover, functions $\Re_{\p P}$ and $\Im_{\p P}$ are non-constant, bounded, Dirichlet-finite, and harmonic for the random walk $(\compk, p_{\p P})$.
\endproclaim

\heading 4.~Random walks of non-univalent circle packings \endheading

We will generalize the notion of random walks induced by \cp s to non-univalent \cp s.
For branched packings without angles of overlaps (see [D1], [D2], [DSt]) one can use Definition~3.1.
However, in order to have a comprehensive approach, we proceed as follows.

Suppose \mmp P is a packing for \mcompk, possibly branched and with prescribed angles of overlaps (see [BoSt], [T2]).
Let $v\in \intk$.
As before, if $u$ is a neighbor of $v$ then there are exactly two vertices $w'$ and $w''$ that are adjacent to both $u$ and $v$.
Let $z'$ and $z''$ be the radical centers (see [C], [Y]) of triples of circles $C_{\p P}(u)$, $C_{\p P}(v)$, $C_{\p P}(w')$ and $C_{\p P}(u)$, $C_{\p P}(v)$, $C_{\p P}(w'')$, respectively.

\midinsert
\epsfysize=8.5truecm
\centerline{\epsffile[210 300 400 490]{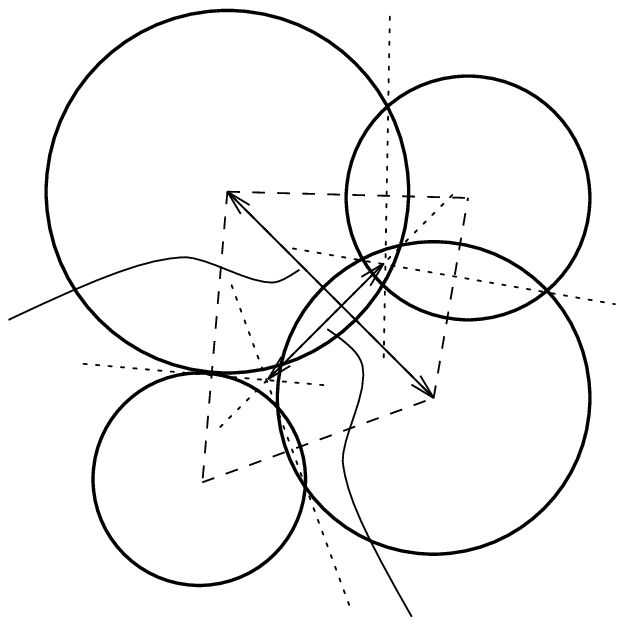}}
\vskip-8.0truecm\hskip5.0truecm $\circle Pu$ 
\vskip0.0truecm\hskip10.5truecm $C_{\p P}(w'')$
\vskip2.4truecm\hskip0.2truecm $|\smap P(u) - \smap P(u)|$
\vskip2.1truecm\hskip3.1truecm $C_{\p P}(w')$ \hskip3.7truecm $\circle Pv$
\vskip0.7truecm\hskip8.6truecm $|z'-z''|$
\captionwidth{25pc}
\botcaption{Figure~3} Radical centers and the edge conductance.
\endcaption
\endinsert

We set the edge conductance between $v$ and $u$ induced by \mmp P to be
$$
c_{\p P}(u,v) := \frac{|z'-z''|}{|\smap P(u) - \smap P(u)|}.
$$
We are now ready to describe the random walk of \mmp P.

\definition{Definition~4.1}
If \mmp P is a \cp\ for \mcompk, possibly branched and with prescribed angles of overlaps, then the random walk $(\compk, p_{\p P})$ induced by \mmp P is defined as follows:
$$
p_{\p P}(u,v)=\cases \frac{c_{\p P}(u,v)}{\sum_{w\sim u}c_{\p P}(u,w)}, &\text{if $v\in \intk$, $u\sim v$, and $u\ne v$}\\
	0, &\text{if $v\in \intk$ and $u=v$, or $v\notin \intk$ and $u\ne v$}\\
	1, &\text{if $v\notin \intk$ and $u=v$}.\endcases
$$
\enddefinition

\remark{Remark~4.2}
If \mcompk\ has no boundary then it follows from the definition that $(\compk, p_{\p P})$ is reversible with edge conductances $c_{\p P}(u,v)$.
If all angles of overlaps in \mmp P are 0 (i.e., \mmp P is a locally univalent or branched \cp) then $c_{\p P}(u,v)$ are given by the formula in Remark~3.2.2).
\endremark

We finish with the following generalization of Lemma~3.3; the proof is similar to that of Lemma~3.3 and is left to the reader.

\proclaim{Lemma~4.3}
Suppose \mmp P is a \cp\ for \mcompk, possibly branched and with prescribed angles of overlaps.
Then both functions $\Re_{\p P}, \Im_{\p P}: \vertk \to \bold R$, $\Re_{\p P}(v):=\text{{\rm Re}}\bigl(\smap P(v)\bigr)$, $\Im_{\p P}(v):=\text{{\rm Im}}\bigl(\smap P(v)\bigr)$, are harmonic for $(\compk, p_{\p P})$.
\endproclaim

\enddocument

%% file: randomw_defs.tex
\define\cp{circle packing}

\define\compk{\Bbb K}
\define\mcompk{$\compk$}
\define\vertk{{\Bbb K}^0}
\define\mvertk{${\Bbb K}^0$}
\define\edgek{{\Bbb K}^1}
\define\medgek{${\Bbb K}^1$}
\define\intk{\text{\rm{int}}\, {\Bbb K}^0}
\define\mintk{$\intk$}

\define\p#1{\Cal #1}
\define\mmp#1{$\p #1$}

\define\circle#1#2{C_{\p #1}(#2)}
\define\mcircle#1#2{$\circle #1#2$}

\define\r#1#2{r_{\p #1}(#2)}
\define\mr#1#2{$\r #1#2$}


\define\smap#1{S_{\p #1}}
\define\msmap#1{$\smap #1$}

\define\carr#1{\text{\rm{carr}} (\p #1)}
\define\mcarr#1{$\carr #1$}